\documentclass[11pt,twoside]{article}
\usepackage{amsfonts}
\usepackage{mathrsfs}
\usepackage{bbm}
\usepackage{amsfonts}
\usepackage{amssymb,amsmath,graphicx}
\usepackage{float}
\usepackage[colorlinks=true]{hyperref}

\hypersetup{urlcolor=blue, citecolor=red}

\allowdisplaybreaks

\begin{document}
\title{ Protection zone in a diffusive predator-prey model with Beddington-DeAngelis functional
response\thanks{Supported by the National Natural Science Foundation
of China (11171048).}}
\author{Xiao He$^{\rm a,b}$ \quad Sining Zheng$^{\rm b,}${\thanks{Corresponding author.
E-mail: xiaohemath@163.com\,(X. He), snzheng@dlut.edu.cn\,(S.N. Zheng)}}\\
\footnotesize $^{\rm a}$Department of Mathematics, Dalian Minzu University, Dalian 116600, P.R. China\\
\footnotesize$^{\rm b}$School of Mathematical Sciences, Dalian University of
Technology, Dalian 116024, P. R. China} \maketitle
\date{}
\newtheorem{theorem}{Theorem}
\newtheorem{definition}{Definition}[section]
\newtheorem{lemma}{Lemma}[section]
\newtheorem{proposition}{Proposition}[section]
\newtheorem{corollary}{Corollary}[section]
\newtheorem{remark}{Remark}
\renewcommand{\theequation}{\thesection.\arabic{equation}}
\catcode`@=11 \@addtoreset{equation}{section} \catcode`@=12
\maketitle{} \begin{abstract} In any reaction-diffusion system of
predator-prey models, the population densities of species are
determined by the interactions between them, together with the
influences from the spatial environments surrounding them.
Generally, the prey species would die out when their birth rate
is too low, the habitat size is too small, the predator grows
too fast, or the predation pressure is too high. To save the endangered
prey species, some human interference is useful, such as creating a
protection zone where the prey could cross the boundary freely but
the predator is prohibited from entering. This paper studies the existence of
positive steady states to a predator-prey model with reaction-diffusion terms, Beddington-DeAngelis type functional response and non-flux boundary conditions. It is shown that there is a
threshold value $\theta_0$ which characterizes the refuge ability of prey such that the positivity of prey
population can be ensured if either the prey's birth  rate satisfies $\theta\geq\theta_0$ (no matter how large the predator's growth rate is) or the predator's growth rate satisfies
$\mu\le 0$, while a protection zone $\Omega_0$ is necessary for such
positive solutions if $\theta<\theta_0$ with $\mu>0$ properly large. The more
interesting finding is that there is another threshold value
$\theta^*=\theta^*(\mu,\Omega_0)<\theta_0$, such that the positive solutions do exist for all
$\theta\in(\theta^*,\theta_0)$. Letting $\mu\rightarrow\infty$, we get the third threshold value $\theta_1=\theta_1(\Omega_0)$  such that if $\theta>\theta_1(\Omega_0)$, prey species could survive no matter how large the predator's growth rate is. In addition, we get the fourth threshold value $\theta_*$ for negative $\mu$ such that the system
admits positive steady states if and only if $\theta>\theta_*$.
 All these results match well with the mechanistic derivation for the B-D type
functional response recently given by Geritz and Gyllenberg [A
mechanistic derivation of the DeAngelis-Beddington functional
response, J. Theoret. Biol. 314 (2012) 106--108]. Finally, we obtain the uniqueness of positive steady states for $\mu$ properly large,
as well as the asymptotic behavior of the unique positive steady
state as $\mu\rightarrow\infty$.
\begin{description}
\item[MSC:]  92D40, 35J47, 35K57
\item[Keywords:] Reaction-Diffusion; Predator-Prey; Beddington-DeAngelis type functional
response; Protection zone; Bifurcation
\end{description}
\end{abstract}


\section{Introduction}

\mbox\indent Biological resources are renewable, but many
have been exploited unreasonably. Nowadays, some species cannot survive in their habitat without human intervention. Such interventions have included establishing banned fishing areas and fishing periods to cope
with over-fishing in fishery production, and setting up nature reserves to protect endangered species. These phenomena are usually described via diffusive predator-prey models, where the population evolution of the
species relies on the interactions between predator and prey, as well as the
influences from the spatial environments surrounding them.
Naturally, prey species would die out when the prey's birth rate is too
low, the habitat size is too small, the predator's growth rate is too
fast, or the predation rate is too high. To save the endangered prey
species, various human interferences are proposed such as creating a
protection zone where the prey could cross the boundary freely but
the predator is prohibited from entering. Refer to the works on
protection zones by Du {\it et al} for the Lotka-Voltera type
competition system \cite{DL}, Holling II type predator-prey system
\cite{DS2006}, Leslie type predator-prey system \cite{DPW}, as well
as predator-prey systems with protection coefficients \cite{DS2007}.
Oeda studied the effects of a cross-diffusive Lotka-Voltera type
predator-prey system  with a protection zone \cite{O}. A
cross-diffusive Lotka-Voltera type competition system with a
protection zone was investigated by Wang and Li \cite{WL}. Zou and
Wang studied an ODE model of protection zones, where the sizes of the protection zones are reflected by
restricting the functionals' coefficient for the predator \cite{ZW}.
Recently, Cui, Shi and Wu observed the  strong Allee effect in a diffusive predator-prey system with protection zones \cite{RC}.

In this paper, we study the steady states to the following diffusive
predator-prey system with Beddington-DeAngelis type functional response
\begin{eqnarray}\label{a} \left\{
\begin{array}{llll}\displaystyle u_t-d_1\Delta
u=u(\theta-u-\frac{a(x)v}{1+mu+kv}),&x\in \Omega,~~t>0,
\\[6pt]
\displaystyle v_t-d_2\Delta v=v(\mu-v+\frac{cu}{1+mu+kv}),& x\in \Omega_1,~~t>0,\\[4pt]
\displaystyle \frac{\partial u}{\partial {n}}=0 ,& x\in\partial\Omega,~~t>0,\\[4pt]
\displaystyle\frac{\partial v}{\partial {n}}=0 ,& x\in\partial\Omega_1,~~t>0,\\[4pt]
\displaystyle u(x,0)=u_0(x)\geq(\not\equiv) 0,&x\in\Omega,\\[4pt]
\displaystyle v(x,0)=v_0(x)\geq(\not\equiv) 0, & x\in \Omega_1,
\end{array}\right. \end{eqnarray}  
where $\Omega$ is a bounded domain in $\mathbb{R}^N$ $(N\leq3)$ with
smooth boundary $\partial\Omega$, ${\Omega}_0\Subset\Omega$
with $\partial \Omega_0$ smooth,
$\Omega_1=\Omega\backslash\overline{\Omega}_0$, constants
$d_1,d_2,\theta,c,m,k>0$, $\mu\in \mathbb{R}$,
$\frac{\partial}{\partial n}$ is the outward normal derivative on the
boundary, and\begin{eqnarray} a(x)=\left\{
\begin{array}{llll}
0,&x\in\overline{\Omega}_0,\\a,&x\in\Omega_1.\end{array}\right.\end{eqnarray}
The fact of $a(x)=0$ in ${\Omega}_0$ implies that no predation could take place there.

 Eq.\,(\ref{a}) is a reaction-diffusion system of species $u$ and $v$,
and the dynamical behavior of species would be determined not only
by the mechanism of the functional response between $u$ and $v$, but
also by the interaction between their reaction and diffusion. Here
prey $u$ and predator $v$ disperse at rates $d_1$ and $d_2$, and
grow at rates $\theta$ and $\mu$, respectively. The prey is consumed
with the functional response of Beddington-DeAngelis type
$\frac{a(x)uv}{1+mu+kv}$ in $\Omega$, and contributes to the
predator with growth rate $\frac{cuv}{1+mu+kv}$ in $\Omega_1$.
Non-flux boundary conditions mean that the habitat of the two species is
closed. The B-D type functional response was introduced by
Beddington \cite{B} and DeAngelis \cite{DGO}. Refer to
\cite{B,DGO,DK} for the background of the original predator-prey
model with B-D type functional response. Guo and Wu studied the
existence, multiplicity, uniqueness and stability of the positive
solutions under homogeneous Dirichlet boundary conditions in
\cite{GW2010}, as well as the effect of large $k$ in \cite{GW2012}.
Chen and Wang established the existence of nonconstant positive
steady-states under Neumann boundary conditions \cite{CW,PW}.

In particular, a mechanistic derivation for the B-D type functional
response has been given by Geritz and Gyllenberg in
\cite{GG} recently, where predators $v$ were divided into
searchers $v_S$ with attack rate $a$ and handlers $v_H$ with
handling time $h$, while preys $u$ were structured
into two classes: active preys $u_P$ and those prey individuals
$u_R$ who have found a refuge with total refuge number $b$ and
 sojourn time $\tau$. In these terms, the parameters in B-D
type functional response of (\ref{a}) can be understood as that
$m=ah$ reflects the handling time of $v_H$, and $k=b\tau$
describes the refuge ability of the prey.

The prey's refuge may come from its aggregation, reduction of its activity, or places where its
predation risk is somehow reduced \cite{S}. Dynamic consequences of prey refuges were observed by
Gonz\'{a}lez-Olivares and Ramos-jiliberto with more prey, fewer predators and enhanced stability \cite{GR}.
On the other hand, refuges from species usually cost the prey in terms of reduced
feeding or mating opportunities \cite{S}, and hence their population could
not be very large. In contrast, the protection zones, as refuges from humans, always benefit the endangered species. Refer to \cite{HRVL,KR,MD,SMR,WW, WF, YR} for more backgrounds on prey refuges and their affections. In this paper, we will show the effect of the prey's refuge and the size of the protection zone on the coexistence and stability of the predator-prey system with B-D type functional response.
The results obtained here observe the general law that refuges and protection zones benefit the coexistence of species \cite{GR,S,ZW}.

Since the model (\ref{a}) contains different coefficients $a(x)$ and $c$ in the B-D type functional response terms for $u$ and $v$ respectively, without loss of generality, suppose $d_1=d_2=1$ for simplicity. The steady-state problem corresponding to (\ref{a}) takes the form 
\begin{eqnarray} \left\{ \begin{array}{llll}
\displaystyle-\Delta u=u(\theta-u-\frac{a(x)v}{1+mu+kv})&\mbox{in}~\Omega,\\[6pt]
\displaystyle-\Delta v=v(\mu-v+\frac{cu}{1+mu+kv})&\mbox{in}~\Omega_1,\\[4pt]
\displaystyle\frac{\partial u}{\partial
{n}}\Big|_{\partial\Omega}=0,~~ \frac{\partial v}{\partial
{n}}\Big|_{\partial\Omega_1}=0.&\end{array}\right.
\label{s}\end{eqnarray} 


Denote by $\lambda_1(q)$ the first eigenvalue of $-\Delta+q$ over
$\Omega$ under homogeneous Neumann boundary conditions with
$q=q(x)\in L^\infty(\Omega)$. The following properties of
$\lambda_1(q)$ are well known:
\begin{itemize}
\item[(i)] $\lambda_1(0)=0$;
\item[(ii)] $\lambda_1(q_1)>\lambda_1(q_2)$ if $q_1\geq q_2$ and $q_1\not\equiv q_2$;
\item[(iii)] $\lambda_1(q)$ is continuous with respect to $q\in L^\infty(\Omega)$.
\end{itemize}

Define
\begin{align} \label{q2}
\theta^*(\mu,\Omega_0)=\lambda_1(q(x)),~~\theta_0=\frac{a}{k},~~\theta_1(\Omega_0)=\lambda_1(q_0(x)),
\end{align}
with
\begin{align}\label{q1}q(x)=\frac{a(x)\mu}{1+k\mu},~~q_0(x)=
\left\{\begin{array}{lll}0, &x\in \overline\Omega_0,\\
\theta_0, &x\in\Omega_1.
\end{array}\right.\end{align}
Denote by $U_{\theta,q_0}$ the solution of the scalar problem
\begin{align}
\label{uq}-\Delta u=u(\theta-u-q_0(x))~{\rm
in}~\Omega,~~\frac{\partial u}{\partial{n}}=0~{\rm
on}~\partial\Omega.
\end{align}

Due to $\theta_1=\inf\limits_{\phi\in
H^1(\Omega),\int_{\Omega}\phi^2dx>0}\frac{\int_{\Omega}|\nabla\phi|^2dx+\frac{a}{k}\int_{\Omega_1}\phi^2dx}{\int_{\Omega}\phi^2dx}$,
the properties (i)--(iii) of $\lambda_1(q)$ imply the following lemma:
\begin{lemma}\label{l1.1}
 $\theta^*{(\mu,\Omega_0)}$ is strictly increasing with respect to
$\mu$ and decreasing when $\Omega_0$ enlarging, $\theta^*(0,\Omega_0)=0$, $\theta^*(\mu,\Omega_0)<\theta_0$,
$\lim_{\mu\rightarrow\infty}\theta^*(\mu,\Omega_0)=\theta_1(\Omega_0)\leq\frac{a|\Omega_1|}{k|\Omega|}$,
$\lim_{|\Omega_1|\rightarrow0}\theta^*(\mu,\Omega_0)=0$,
$\lim_{|\Omega_0|\rightarrow0}\theta^*(\mu,\Omega_0)=\frac{a\mu}{1+k\mu}$.\qquad$\Box$
\end{lemma}

Biologically, we are interested in the positivity of the prey $u$ in
the diffusive predator-prey model (\ref{s}). We state the main
results of the paper one by one as follows.\medskip

 Obviously, either large $\theta$ or small $\mu$
benefits the prey $u$. In the first theorem, we give two sufficient
conditions for keeping the prey positive without
protection zones.
\begin{theorem}\label{th1}
If $\theta\ge \theta_0$ or $\mu\le0$, then the positivity of $u$
would be ensured automatically without any protections zones.
\end{theorem}

The next theorem implies that a suitable protection zone guarantees the
existence of positive solutions to (\ref{s}) under $\theta<\theta_0$ with $\mu>0$.
\begin{theorem}\label{th2}
Suppose $\mu> 0$. If $\theta^*(\mu,\Omega_0)< \theta< \theta_0$,
then Eq.\,(\ref{s}) has at least one positive solution. Furthermore,
if $\theta\leq\theta^*(\mu,\Omega_0)$ with $m\le\frac{(k\mu+1)^2}{a\mu}$,
then Eq.\,(\ref{s}) has no positive solutions.
\end{theorem}

 In the third theorem, we give a necessary and sufficient condition for the coexistence of $u$ and $v$ under $\mu\in(-\frac{c}{m},0]$.

\begin{theorem}\label{th2'} Suppose $-\frac{c}{m}<\mu\le 0$. Then
Eq.\,(\ref{s}) has at least one positive solution if and only if
$\theta>\theta_*=-\frac{\mu}{c+m\mu}=\frac{|\mu|}{c-m|\mu|}\ge
0$.\end{theorem}

\begin{remark}\label{mk1}{\rm Since $\lim_{|\Omega_1|\rightarrow0}\theta^*(\mu,\Omega_0)=0$ by Lemma \ref{l1.1},
for any $\theta>0$ and $\mu\ge 0$ fixed, the key condition
$\theta>\theta^*(\mu,\Omega_0)$ in Theorem \ref{th2}  can be
realized by enlarging the size of the protection zone
$\Omega_0=\Omega\setminus\overline \Omega_1$. So does the condition
$\theta>\theta_1$ in the following Theorem \ref{th3}. }\end{remark}

\begin{remark}\label{mk2}
{\rm Theorem \ref{th1} shows that no protection zones are necessary
for the positivity of $u$ if $\mu\le 0$. It is
known by Theorem \ref{th2'} that in addition to the positivity of
$u$, the positivity of $v$ can be ensured also if the death rate of
the predator $v$ is not too high with $\mu\in (-\frac{c}{m},0]\subset
(-\infty,0]$ and the birth rate of the prey $u$ is properly large such that
$\theta>\theta_*$.}
\end{remark}

Finally, the last theorem says the positive solutions of (\ref{s})
are in fact unique if $\theta$ is even larger than $\theta_1$
under large $\mu$, and determines the asymptotic behavior of the
unique positive solution as $\mu\rightarrow \infty$. In fact, from Lemma \ref{l1.1} and Theorem \ref{th2} that
if $\theta>\theta_1$,  prey species could be alive  no matter how large the predator's growth rate is.

\begin{theorem}\label{th3}
 If $\theta>\theta_1(\Omega_0)$, then there exists $\mu^*>0$ such that the positive solution of (\ref{s}) is
 unique and linearly stable when $\mu\geq\mu^*$. Furthermore, the
 unique positive solution satisfies
$(u,v-\mu)\rightarrow (U_{\theta,q_0},0)$ uniformly on
$\overline\Omega$ and $\overline\Omega_1$, respectively, as
$\mu\rightarrow\infty$.
\end{theorem}

This paper is arranged as follows. In the next two sections, we
give the proofs of Theorems \ref{th1}--\ref{th2'} and Theorem \ref{th3},
respectively. The last section is devoted to a discussion of the obtained results, by analyzing them with the mechanistic derivation for the B-D type functional response in \cite{GG}.

\section{Existence of positive solutions}
\mbox\indent At first we deal with the proof of Theorem
\ref{th1}.\medskip

\noindent{\bf Proof of Theorem \ref{th1}.}  Assume
$\mu\leq-\frac{c}{m}$. Integrate the second equation of (\ref{s})
over $\Omega_1$ to get
$$0=\int_{\Omega_1}v\Big(\mu-v+\frac{cu}{1+mu+kv}\Big)dx,$$and hence
$$0\leq\int_{\Omega_1}v^2dx=\int_{\Omega_1}v(\mu+\frac{cu}{1+mu+kv})dx\leq(\mu+\frac{c}{m})\int_{\Omega_1}vdx\leq
0.$$This concludes $v\equiv0$, and so $u$
satisfies\begin{align}\label{uu} -\Delta
u=u(\theta-u)~~\mbox{in}~\Omega,~~~\frac{\partial u}{\partial
n}=0~~{\rm on}~\partial\Omega.\end{align} Obviously, (\ref{uu})
admits the solution $u=\theta>0$.\medskip

The desired result for $-\frac{c}{m}<\mu\le 0$ is substantially
concluded from Theorem \ref{th2'}. Indeed, the subcase of
$\theta>\theta_*$ is covered by Theorem \ref{th2'}, while for
$\theta\le \theta_*$, it can be found in the proof of Theorem
\ref{th2'} that $v\equiv 0$, and so $u=\theta>0$.\medskip

Next consider the first equation of (\ref{s}) with
$\theta\geq\theta_0$. It is easy to know that
$\frac{a(x)v}{1+mu+kv}<\frac{a(x)}{k}\le
\frac{a}{k}=\theta_0$ for $v\ge 0$. Thus, for any $v(x)\ge 0$,
there is $\tilde\theta_0\in (0,\theta_0)$ such that
\begin{align*}\label{uuu}
-\Delta u=u(\theta-u-\frac{a(x)v}{1+mu+kv})>
u(\theta-\tilde\theta_0-u)~~\mbox{in}~\Omega,~~\frac{\partial
u}{\partial n}=0~~{\rm on}~\partial \Omega.
\end{align*}
This ensures that $u\ge
\theta-\tilde\theta_0>0$.\qquad$\Box$\medskip

We need some preliminaries represented as lemmas and propositions for the proof of Theorem \ref{th2},
and begin with two known results on the maximum principle and the
Harnack inequality.
\begin{lemma}\label{l2.1} (Maximum Principle \cite{LN}) Let
$g\in C(\overline{\Omega}\times \mathbb{R})$, $w\in
C^2(\Omega)\bigcap C^1(\overline\Omega)$, where $\Omega$ is a
bounded domain in $\mathbb{R}^N$ with smooth boundary.
\begin{description}
\item{ {\rm (a)}} If ~$\Delta w+g(x,w)\leq0 $ in
$\Omega$, $\frac{\partial w}{\partial {n}}\geq0$ on  $\partial
\Omega$ and $\min_{\overline \Omega} w=w(x_0)$, then
$g(x_0,w(x_0))\leq0$.
\item{ {\rm (b)}} If ~$\Delta w+g(x,w)\geq0$ in $\Omega$, $\frac{\partial
w}{\partial {n}}\leq0$ on $\partial \Omega$ and $\max_{\overline
\Omega} w=w(x_0)$, then $g(x_0,w(x_0))\geq0$.
\end{description} \quad$\Box$
\end{lemma}
\begin{lemma}\label{l2.2} (Harnack Inequality
\cite{LN1999,LNT}) Let $f\in L^p(\Omega)$ with
$p>\max\{\frac{N}{2},1\}$, and w be a non-negative solution of
$\Delta w+f(x)w=0$ in a bounded domain $\Omega\subset\mathbb{R}^N$ with smooth
boundary  under homogeneous Neumann boundary condition. Then there exists a positive constant $C=C(p,N,\Omega,\|
f\|_{L^p(\Omega)})$ such that
$$\max\limits_{\overline\Omega}w\leq C\min\limits_{\overline\Omega}w.\qquad\Box$$
\end{lemma}

The following {\it a priori} estimates are easy to get.
\begin{lemma}\label{l2.3}  Let $(u,v)$ be a nontrivial non-negative solution
of (\ref{s}). Then
\begin{eqnarray*}
0<u\leq\theta,~~\mu_+<v\leq\mu_++\frac{c\theta}{1+m\theta+k\mu_+},~~
\|u\|_{C^{1,\alpha}(\overline\Omega)}+\|v\|_{C^{1,\alpha}(\overline\Omega_1)}\leq
C,\end{eqnarray*} with $\mu_+=\max\{\mu,0\}$, $\alpha\in(0,1)$,
$C=C(\theta,\mu,\Omega_0)>0$.
\end{lemma}
{\bf Proof.} Suppose $u(x_0)=\max_{\overline \Omega}u(x)>0$. By
Lemma \ref{l2.1}(b), we have
$$u(x_0)\Big(\theta-u(x_0)-\frac{a(x_0)v(x_0)}{1+mu(x_0)+kv(x_0)}\Big)\geq
0,$$ and then
$$ u(x_0)\leq\theta-\frac{a(x_0)v(x_0)}{1+mu(x_0)+kv(x_0)}\leq\theta.$$  Due to Lemma \ref{l2.2}, we arrive at
$0<u\leq\theta$ on $\overline\Omega$. Similarly, we can show
$\mu_+<v\leq\mu_++\frac{c\theta}{1+m\theta+k\mu_+}$ on
${\overline\Omega}_1$.

 The $C^{1,\alpha}$ boundedness of  solutions comes from
the elliptic regularity theory together with the Sobolev embedding
theorem.\qquad$\Box$\medskip

We will use the local bifurcation theorem of Crandall and Rabinowitz
\cite{CR} and the global bifurcation theorem of Rabinowitz \cite{R}
to prove Theorem \ref{th2}.

Denote the semitrivial solution curves by
\begin{eqnarray*}
\Gamma_u=\{(\theta,u,v)=(\theta,0,\mu);\,\theta>0\},~~
\Gamma_v=\{(\theta,u,v)=(\theta,\theta,0);\,\theta>0\}.
\end{eqnarray*}
Define \begin{eqnarray*} X=W^{2,p}_n(\Omega)\times
W^{2,p}_n(\Omega_1),~Y=L^p(\Omega)\times L^p(\Omega_1)~{\rm with}~p>N,~Z=C_n^1(\overline\Omega)\times
C_n^1(\overline\Omega_1),
\end{eqnarray*}
where \begin{eqnarray*} W^{2,p}_n(\Omega)=\{w\in
W^{2,p}(\Omega);~\frac{\partial w}{\partial n}=0~{\rm on}~\partial
\Omega\},~C_n^1(\overline\Omega)=\{w\in
C^1(\overline\Omega);\frac{\partial w}{\partial n}=0~{\rm
on}~\partial\Omega\}.
\end{eqnarray*}
The Sobolev embedding theorem implies $X\subseteq Z$.

Let $(\phi^*,\psi^*)$ solve
\begin{align*}
&\Delta\phi^*+(\theta^*-\frac{a(x)\mu}{1+k\mu})\phi^*=0~{\rm
in}~\Omega,~\frac{\partial \phi^*}{\partial
 {n}}=0~{\rm on}~\partial\Omega,\\
&\Delta\psi^*-\mu\psi^*+\frac{c\mu}{1+k\mu}\phi^*=0~{\rm
in}~\Omega_1,~\frac{\partial \psi^*}{\partial
 {n}}=0~{\rm on}~\partial\Omega_1.
\end{align*}
Then $\psi^*=(-\Delta+\mu I)^{-1}_{\Omega_1}\frac{c\mu}{1+k\mu}\phi^*$.
\begin{proposition}\label{p2.1} Let $\mu>0$. Then there are positive solutions
of (\ref{s}) bifurcating from $\Gamma_u$ if and only if
$\theta>\theta^*(\mu,\Omega_0)$, possessing the form
\begin{eqnarray}\
\Gamma_1=\{(\theta,u,v)=(\theta(s),s\phi^*+o(|s|),\mu+s\psi^*+o(|s|));\,s\in
 (0,\sigma)\}\end{eqnarray} with
 $(\theta(0),u(0),v(0))=(\theta^*,0,\mu)$ for some $\sigma>0$ in a neighborhood
 of  $(\theta^*,0,\mu)\in \mathbb{R}\times X$.
\end{proposition}
{\bf Proof.}  Denote by $V=v-\mu$,
\begin{eqnarray}
F(\theta,u,V)=\Big(\begin{array}{clcr} \ \Delta u+f_1(\theta,u,V+\mu)\\
\Delta V+f_2(\mu,u,V+\mu)\end{array}\Big)^T~{\rm and}~
F_1(\theta,u,v)=\Big(\begin{array}{clcr} \ \Delta u+f_1(\theta,u,v)
\\
\Delta v+f_2(\mu,u,v)
\end{array}\Big)^T\end{eqnarray}
with
\begin{align*}f_1(\theta,u,v)=u(\theta-u-\frac{a(x)v}{1+mu+kv}),~~
f_2(\mu,u,v)=v(\mu-v+\frac{cu}{1+mu+kv}).\end{align*} Obviously,
 $F(\theta,u,V)=0 $ is equivalent to $ F_1(\theta,u,v)=0$, and
$F_1(\theta,0,\mu)=F(\theta,0,0)=0$ for $\theta\in \mathbb{R}$.
A direct calculation yields
\begin{eqnarray}\
F_{(u,V)}(\theta,0,0)[\phi,\psi]=\Big(\begin{array}{clcr}\
\Delta\phi+(\theta-\frac{a(x)\mu}{1+k\mu})\phi
\\ \Delta\psi-\mu\psi+\frac{c\mu}{1+k\mu}\phi
\end{array}\Big)^T.\end{eqnarray}\
By the Krein-Rutman theorem, $F_{(u,V)}(\theta,0,0)[\phi,\psi]=(0,0)$
has a solution  $\phi>0$ if and only if $\theta=\theta^*$. So
$(\theta^*,0,\mu)$ is the only possible bifurcation point from which
positive solutions of (\ref{s}) bifurcate from $\Gamma_u$. Besides, we have
\begin{align*}
 {\rm Ker}
 F_{(u,V)}(\theta^*,0,0)={\rm Span}\,\{(\phi^*,\psi^*)\},~~
 \makebox{dim\,Ker}\,
 F_{(u,V)}(\theta^*,0,0)=1.
 \end{align*}
 For $(\bar\phi,\bar\psi)\in Y\cap
 \makebox{Range}\,F_{(u,V)}(\theta^*,0,0)$,  choose $(\phi,\psi)\in
 X$ such that
\begin{eqnarray}\label{ee1} \left\{ \begin{array}{llll}
\displaystyle\Delta\phi+(\theta-\frac{a(x)\mu}{1+k\mu})\phi=\bar\phi,\\[6pt]
\displaystyle\Delta\psi-\mu\psi+\frac{c\mu}{1+k\mu}\phi=\bar\psi.
\end{array}
\right.\end{eqnarray} Multiplying by $\phi^*$ on both sides of the first equation of
(\ref{ee1}) and integrating by parts over $\Omega$, we get
$\int_\Omega\bar\phi\phi^* dx=0$. Then
\begin{eqnarray}
\makebox{Range}\,F_{(u,V)}(\theta^*,0,0)=\Big\{(\bar\phi,\bar\psi)\in
Y;\int_\Omega\bar\phi\phi^* dx=0\Big\},
\end{eqnarray}
and thus
$$\makebox{codim\,Range}\,F_{(u,V)}(\theta^*,0,0)=1.$$
By a simple calculation,
\begin{align*}
&F_\theta(\theta^*,0,0)=F_{\theta\theta}(\theta^*,0,0)=(0,0),\\&F_{\theta(u,V)}(\theta^*,0,0)[\phi^*,\psi^*]=
(\phi^*,0) \notin
\makebox{Range\,}F_{(u,V)}(\theta^*,0,0).\end{align*}
In conclusion, the proposition  is proved by the local bifurcation theorem \cite{CR}.
\qquad$\Box$
\begin{proposition}\label{p2.2} Let $-\frac{c}{m}<\mu<0$. Then there are positive
solutions of (\ref{s}) bifurcating from $\Gamma_v$ if and only if
$\theta>\theta_*=-\frac{\mu}{c+m\mu}$, having the form
\begin{eqnarray}
\Gamma_2=\{(\theta,u,v)=(\tilde{\theta}(s),\theta+s\phi_*(x)+o(|s|),s+o(|s|));s\in(0,\tilde{\sigma})\}
\end{eqnarray}
 with $\tilde{\theta}(0)=-\frac{\mu}{c+m\mu}$, $\phi_*=(\Delta-\theta I)^{-1}\frac{a(x)\theta}{1+m\theta}$ for some $\tilde{\sigma}>0$ near $(\theta,\theta,0)\in \mathbb{R}\times X$.
\end{proposition}
{\bf Proof.} Let $w=u-\theta$, \begin{eqnarray}\
G(\theta,w,v)=\Big(\begin{array}{clcr}\ \Delta
w+(w+\theta)(-w-\frac{a(x)v}{1+m(w+\theta)+kv})\\
\Delta v+v(\mu-v+\frac{c(w+\theta)}{1+m(w+\theta)+kv})
\end{array}\Big)^T,
\end{eqnarray}\
Then $F_1(\theta,u,v)=0$ is equivalent to $G(\theta,w,v)=0$. We
have
\begin{align*}
&G_{(w,v)}(\theta,w,v)[\phi,\psi]=\\
&\left(\begin{array}{clcr}\
\Delta\phi-(2w+\theta)\phi-\frac{a(x)v}{1+m(w+\theta)+kv}\phi+\frac{a(x)mv(w+\theta)}{(1+m(w+\theta)+kv)^2}\phi-\frac{a(x)(w+\theta)(1+m(w+\theta))}{(1+m(w+\theta)+kv)^2}\psi\\[8pt]
\Delta\psi+(\mu-2v)\psi+\frac{c(w+\theta)}{1+m(w+\theta)+kv}\psi-\frac{ckv(w+\theta)}{(1+m(w+\theta)+kv)^2}\psi+\frac{cv(1+kv)}{(1+m(w+\theta)+kv)^2}\phi
\end{array}\right)^T,
\\[6pt]
&G_\theta(\theta,w,v)=
\Big(\begin{array}{clcr}
-w-\frac{a(x)v}{1+m(w+\theta)+kv}+\frac{a(x)mv(w+\theta)}{(1+m(w+\theta)+kv)^2}\\ \frac{cv(1+kv)}{(1+m(w+\theta)+kv)^2}
\end{array}\Big)^T.\end{align*}
The equation $G_{(w,v)}(\theta,0,0)[\phi,\psi]=(0,0)$ is equivalent to 
\begin{eqnarray}\label{ee2} \left\{
\begin{array}{llll}
\Delta\phi-\theta\phi-\frac{a(x)\theta}{1+m\theta}\psi=0~~~~~&{\rm in}~\Omega,\\[3pt]
\Delta\psi+\mu\psi+\frac{c\theta}{1+m\theta}\psi=0~~~~&{\rm in}~\Omega_1,\\[3pt]
\frac{\partial\phi}{\partial{n}}=0~{\rm
on}~\partial\Omega,~~\frac{\partial\psi}{\partial{n}}=0~{\rm
on}~\partial\Omega_1.
\end{array}\right.
\end{eqnarray}  
The second equation of (\ref{ee2}) has a solution $\psi>0$ if and
only if $\mu=-\frac{c\theta}{1+m\theta}$, i.e. $\theta=-\frac{\mu}{c+m\mu}=\theta_*$. Thus
$(\theta_*,\theta,0)$ is the only possible bifurcation point along
$\Gamma_v$, and $\phi_*$ solves the first equation
of (\ref{ee2}) with $\theta=\theta_*$ and $\psi\equiv 1$. It is
easy to verify that
$$\makebox{Ker}\,G_{(w,v)}(\theta_*,0,0)=\mbox{Span\,}\{(\phi_*,1)\},~~
\makebox{dim\,Ker\,}G_{(w,v)}(\theta,0,0)=1.$$A direct calculation
shows
\begin{align*}
&G_\theta(\theta_*,0,0)=G_{\theta\theta}(\theta_*,0,0)=(0,0),\\
&\mbox{Range\,}G_{(w,v)}(\theta_*,0,0)=\{(f,g)\in Y;\int_\Omega
gdx=0\},~\mbox{codim\,Range\,}G_{(w,v)}(\theta_*,0,0)=1,\\
&G_{\theta(w,v)}(\theta_*,0,0)[\phi_*,1]=\big(
-\phi_*-\frac{a(x)}{(1+m\theta)^2},
\frac{c}{(1+m\theta)^2}\big)\notin
\mbox{Range\,}G_{(w,v)}(\theta_*,0,0).
\end{align*}
By the local bifurcation theorem \cite{CR}, we get the desired results of the Proposition \ref{p2.2}.
\qquad$\Box$\medskip

In order to use the global bifurcation theorem for $\mu>0$, define $
F_2:\mathbb{R}\times Z\rightarrow Z$ by
\begin{eqnarray}
\ F_2(\theta,u,v)=
\Big(\begin{array}{clcr}\ u \\
v-\mu\end{array}\ \Big)^T-\Big(\begin{array}{clcr}\ (-\Delta+I)^{-1}_\Omega(u+f_1(\theta,u,v))\\
(-\Delta+I)^{-1}_{\Omega_1}(v-\mu+f_2(\mu,u,v))
\end{array}\ \Big)^T.
\end{eqnarray}
Then (\ref{s}) is equivalent to $F_2(\theta,u,v)=0$. Let
$\tilde{\Gamma}_1\subset \mathbb{R}\times Z$ be the maximal
connected set satisfying \begin{align*}
\Gamma_1\subset\tilde{\Gamma}_1\subset\{(\theta,u,v)\in
\mathbb{R}\times
Z\backslash\{(\theta^*,0,\mu)\};F_2(\theta,u,v)=(0,0)\}.
\end{align*}

From the global bifurcation theory of Rabinowitz \cite{R}, one of
the following non-excluding results must be true (see Theorem 6.4.3
in \cite{G}):
\begin{description}
\item[\rm (a)] $\tilde{\Gamma}_1$
is unbounded in $\mathbb{R}\times Z$.\item[\rm(b)]  There exists a
constant $\bar\theta\neq\theta^*$ such that
$(\bar\theta,0,\mu)\in\tilde \Gamma_1$.\item[\rm (c)]  There exists
$(\tilde{\theta},\tilde{\phi},\tilde{\psi})\in
\mathbb{R}\times(Y_1\backslash\{(0,\mu)\})$ with $Y_1=\{(\bar\phi,\bar\psi)\in
Z;\int_\Omega\bar\phi\phi^*dx=0\}$ such that
$(\tilde{\theta},\tilde{\phi},\tilde{\psi})\in\tilde{\Gamma}_1$.
\end{description}

Now we give the proofs of Theorems \ref{th2} and \ref{th2'}.

\noindent{\bf Proof of Theorem \ref{th2}.} At first we know that
 $u,v>0$ for any $(\theta,u,v)\in\tilde{\Gamma}_1$ which means that the case (c) above cannot occur by $\phi^*>0$.
Otherwise there is a $(\bar\theta,\bar u,\bar v)\in
\tilde{\Gamma}_1$ such that (1) $\bar u>0$ with $\bar v(x_0)=0$ for
some $x_0\in \Omega_1$, or (2) $u(x_1)=v(x_2)=0$ for some $x_1\in
\Omega$ and $x_2\in\Omega_1$, or (3) $\bar v>0$ with $\bar u(x_3)=0$
for some $x_3\in \Omega$. Denote by $\mathscr{B}_\Omega=\{\phi\in
C^1_n(\overline\Omega);\,\phi>0~\mbox{on}~\overline\Omega\}$. Choose
a sequence
$\{(\theta_i,u_i,v_i)\}_{i=1}^\infty\subset\tilde{\Gamma}_1\cap(\mathbb{R}\times\mathscr{B}_\Omega\times\mathscr{B}_{\Omega_1})$
such that
$\lim\limits_{i\rightarrow\infty}(\theta_i,u_i,v_i)=(\bar\theta,\bar
u,\bar v)$ in $\mathbb{R}\times Z$, where $\bar\theta$ can be
$\infty$. Obviously, $(\bar u,\bar v)$ is a non-negative solution of
(\ref{s}) with $\theta=\bar\theta$. By Lemma \ref{l2.2}, one of the
following must hold:
$$\mbox{(1)}~\bar u>0,\bar v\equiv0;~~ \mbox{(2)}~\bar u\equiv0,
\bar v\equiv 0; ~~\mbox{(3)}~ \bar u\equiv0,\bar v>0.
$$

For (3), we have $-\Delta\bar v=\bar v(\mu-\bar v)$ in $\Omega_1$,
$\frac{\partial\bar v}{\partial {n}}=0$ on $\partial\Omega_1$, and
thus $\bar v\equiv\mu$. By Proposition \ref{p2.1}, this implies $\bar\theta=\theta^*$, a contradiction to the definition of
$\tilde{\Gamma}_1$.

Suppose (1) or (2) is true. Integrate the second equation of (\ref{s}) on
$\Omega_1$ with $(u,v)=(u_i,v_i)$ to obtain
\begin{eqnarray*}
\int_{\Omega_1}v_i(\mu-v_i+\frac{cu_i}{1+mu_i+kv_i})dx=0, ~~i\in
\mathbb{N}.
\end{eqnarray*}
On the other hand, $\mu>0$ and $\bar v\equiv0$ ensure $\mu-v_i>0$,
and thus $\mu-v_i+\frac{cu_i}{1+mu_i+kv_i}>0$ for $i$ large enough,
also a contradiction.

The case (b) is excluded by Proposition \ref{p2.1}. So, the only true
case is (a).

From Lemma \ref{l2.3}, $(u,v)$ are uniformly bounded in $Z$ as
$(\theta,u,v)\in\tilde{\Gamma}_1$ which shows that $\theta$ is unbounded.
Combining this with Proposition \ref{p2.1}, we know that (\ref{s})
has at least one positive solution for $\theta>\theta^*(\mu,\Omega_0)$ with $\mu>0$.

Now, let $(u,v)$ be a positive solution of (\ref{s}) with $m\leq\frac{(1+k\mu)^2}{a\mu}$. A direct calculation yields that
$u+\frac{a(x)v}{1+mu+kv}>\frac{a(x)\mu}{1+k\mu}$.
By the monotonicity of the eigenvalue, we conclude that
$$0=\lambda_1\big(-\theta+u+\frac{a(x)v}{1+mu+kv}\big)>\lambda_1\big(-\theta+\frac{a(x)\mu}{1+k\mu}\big).$$
Then
$$\theta>\lambda_1(\frac{a(x)\mu}{1+k\mu})=\theta^*(\mu,\Omega_0).$$
This shows that (\ref{s}) has no positive solution whenever
$\theta\leq\theta^*(\mu,\Omega_0)$ and $m\leq\frac{(1+k\mu)^2}{a\mu}$. \qquad$\Box$ \medskip

\noindent{\bf Proof of Theorem \ref{th2'}}. When $\mu=0$, fix
$\theta>0$. By Lemma \ref{l1.1} and Theorem \ref{th2}, we can take a
sequence $\{(\mu_i,u_i,v_i)\}_{i=1}^\infty$ such that $(u_i,v_i)$ is a
positive solution of (\ref{s}) with $\mu=\mu_i>0$,
$\lim_{i\rightarrow\infty}\mu_i=0$. By Lemma \ref{l2.3} and
embedding theorem, we can choose a subsequence (still denoted by
$\{(\mu_i,u_i,v_i)\}_{i=1}^\infty$) such that $(u_i,v_i)$ converges
to $(\tilde{u},\tilde{v})\in Z$, a non-negative solution of
(\ref{s}). By Lemma \ref{l2.2}, $\tilde{u}>0$ or $\tilde{u}\equiv0$
in $\Omega$; $\tilde{v}>0$ or $\tilde{v}\equiv0$ in $\Omega_1$.

 If $\tilde{u}\equiv0$ and $\tilde{v}>0$, then
$\mu_i-v_i+\frac{cu_i}{1+mu_i+kv_i}<0$ in $\Omega_1$ for $i$ large
enough. This contradicts
$\int_{\Omega_1}v_i(\mu_i-v_i+\frac{cu_i}{1+mu_i+kv_i})dx=0$.

If $\tilde{u}>0$ and $\tilde{v}\equiv0$, then
$\mu_i-v_i+\frac{cu_i}{1+mu_i+kv_i}>0$ in $\Omega_1$ for $i$ large
enough, also a contradiction with
$\int_{\Omega_1}v_i(\mu_i-v_i+\frac{cu_i}{1+mu_i+kv_i})dx=0$.

If $\tilde{u}\equiv0$ and $\tilde{v}\equiv0$, then
$\theta-u_i+\frac{a(x)v_i}{1+mu_i+kv_i}>0$ in $\Omega$ for $i$
large enough,  a contradiction to
$\int_{\Omega}u_i(\theta-u_i+\frac{a(x)v_i}{1+mu_i+kv_i})dx=0$.

In summary, we must have $\tilde{u},\tilde{v}>0$ in $\Omega$
and $\Omega_1$, respectively. This means that (\ref{s}) possesses
positive solutions for all $\theta>0$ if $\mu=0$.

 Now, suppose $-\frac{c}{m}<\mu<0$.  For
$\theta>-\frac{\mu}{c+m\mu}>0$, the existence of positive solutions
can be obtained from Proposition \ref{p2.2} by a similar global bifurcation analysis of $\Gamma_u$ as
that for the branch $\Gamma_v$ with $\mu>0$. We omit the details.

Conversely, let $(u,v)$ be a positive solution of (\ref{s}) with
$\mu\in(-\frac{c}{m},0]$. Then $0<u\leq\theta$ by Lemma \ref{l2.3},
and hence
$$\mu=\lambda_1(v-\frac{cu}{1+mu+kv})>\lambda_1(-\frac{cu}{1+mu})
\geq\lambda_1(-\frac{c\theta}{1+m\theta})=-\frac{c\theta}{1+m\theta},$$
namely, $\theta>-\frac{\mu}{c+m\mu}$. \qquad$\Box$

\section{Uniqueness of positive solutions}
\mbox\indent In this section, we use topological degree to
prove Theorem \ref{th3} for $\theta>\theta_1$ and large $\mu$.
At first, introduce an auxiliary problem
\begin{eqnarray} \left\{ \begin{array}{llll}
\displaystyle-\Delta u=u\Big(\theta-u-\frac{a(x)v}{1+mu+kv}\Big)&\mbox{in}~\Omega,\\[4pt]
\displaystyle-\Delta v=v\Big(\mu-v+t\frac{cu}{1+mu+kv}\Big)&\mbox{in}~\Omega_1,\\[4pt]
\displaystyle\frac{\partial u}{\partial {n}}\Big|_{\partial
\Omega}=0,\quad \frac{\partial v}{\partial {n}}\Big|_{\partial\Omega_1}=0
\end{array}\right. \label{st}\end{eqnarray} 
with  the parameter $t\in[0,1]$. Eq.\,(\ref{st}) reverts back to (\ref{s})
if $t=1$. When $t=0$, we have from the second equation of (\ref{st})
that $v\equiv \mu$, and then obtain the scalar problem
\begin{eqnarray}\label{su} \left\{
\begin{array}{llll}\displaystyle-\Delta u=u(\theta-u-\frac{a(x)\mu}{1+mu+k\mu})&{\rm in}~\Omega,\\
\displaystyle\frac{\partial u}{\partial{n}}=0&{\rm
on}~\partial\Omega,
\end{array}\right.\end{eqnarray}  
which yields Eq.\,(\ref{uq}) as $\mu\rightarrow\infty$.
\begin{lemma}\label{l3.1}
Problem (\ref{uq}) has a unique positive solution if and only if
$\theta>\theta_1$.
\end{lemma}
{\bf Proof.}  Suppose $\theta>\theta_1$.  Let $\phi>0$ be the normalized
eigenfunction with respect to $\theta_1$. Set $\underline{u}=\epsilon\phi$. Then
\begin{align*}
-\Delta\underline{u}=-\epsilon\Delta\phi=\epsilon(\theta_1-q_0(x))\phi
=\epsilon\phi(\theta-q_0(x)-\epsilon\phi)+\epsilon\phi(\theta_1-\theta+\epsilon\phi).
\end{align*}
Choose $\epsilon$ small enough such that
$\theta_1-\theta+\epsilon\phi<0$ to get
$$-\Delta\underline{u}\leq\underline{u}(\theta-q_0(x)-\underline{u})~{\rm in}~\Omega,
~~\frac{\partial\underline{u}}{\partial n}=0~{\rm
on}~\partial\Omega.$$ Obviously, $\underline{u}=\epsilon\phi$ and
$\overline{u}=\theta$ are a pair of positive sub- and supersolutions
of (\ref{uq}) with $\underline{u}\leq\overline{u}$. We can get a positive solution of Eq.\,(\ref{uq}) by the sub-supersolution method. Let $\tilde{u}$
and  $\hat{u}$ be the minimal and maximal positive solutions to
(\ref{uq}), respectively. Since
$$\int_{\Omega}\nabla\tilde{u}\cdot\nabla\hat{u}dx=\int_{\Omega}\tilde{u}\hat{u}(\theta-\tilde{u}-q_0(x))dx=
\int_{\Omega}\tilde{u}\hat{u}(\theta-\hat{u}-q_0(x))dx,$$we conclude
$$\int_{\Omega}\tilde{u}\hat{u}(\tilde{u}-\hat{u})dx=0.$$
Therefore $\tilde{u}\equiv \hat{u}$.

On the other hand, it is obviously true for any positive solution
$u_1$ of (\ref{uq}) that $\theta=\lambda_1(u_1+q_0(x))>\lambda_1(q_0(x))=\theta_1$.
\qquad$\Box$\medskip

Next, we show the uniqueness of positive solutions to (\ref{su}).
\begin{proposition}\label{p3.1}
Suppose $\theta>\theta_1$.  There is a
$\tilde{\mu}=\tilde{\mu}(\theta)>0$ such that for any
$\mu>\tilde{\mu}$, problem (\ref{su}) has an unique positive solution.
\end{proposition}
{\bf Proof.} Since $-q_0(x)<-\frac{a(x)\mu}{1+mU_{\theta,q_0}+k\mu}$, then $U_{\theta,q_0}$ is a subsolution of (\ref{su}). Obviously, $\theta$ is a supersolution of (\ref{su}) and $U_{\theta,q_0}\leq\theta$.  Then there exist
positive solutions to (\ref{su}).

To prove the uniqueness of the positive solutions to (\ref{su}), we at first show that
the positive solutions of (\ref{su}) are linearly stable for large
$\mu$. Let $U$ be a positive solution of (\ref{su}). Consider the
eigenvalue problem
\begin{eqnarray}
-\Delta\phi=\theta\phi-2U\phi-\frac{a(x)\mu(1+k\mu)}{(1+mU+k\mu)^2}\phi+\eta\phi~{\rm
in}~\Omega, ~~\frac{\partial\phi}{\partial{n}}=0~{\rm
on}~\partial\Omega
\end{eqnarray}
with the principal eigenvalues denoted by
\begin{eqnarray}\label{ee} \eta=\eta(\mu)=\inf_{\phi\in H^{1}(\Omega),\,\|
\phi\|_2=1}\int_{\Omega}[|\nabla\phi|^2-\theta\phi^2+2U\phi^2+\frac{a(x)\mu(1+k\mu)}{(1+mU+k\mu)^2}\phi^2]dx.
\end{eqnarray}\
We have
\begin{align*}0=\lambda_1\big(-\theta+2U+\frac{a(x)\mu(1+k\mu)}{(1+mU+k\mu)^2}-\eta\big)>\lambda_1(-\theta-\eta)=-\theta-\eta,\end{align*}
i.e., $\eta>-\theta$. Denote by $\eta^*$ the principal eigenvalue of the problem
\begin{eqnarray}\label{eee}
-\Delta\phi=\theta\phi-2U_{\theta,q_0}(x)\phi-q_0(x)\phi+\eta^*\phi~{\rm
in}~\Omega,~~ \frac{\partial \phi}{\partial{n}}=0~{\rm
on}~\partial\Omega
\end{eqnarray} with the normalized eigenfunction $\phi^*>0$. Then $\eta^*=\frac{\int_\Omega U_{\theta,q_0}^2\phi^* dx}{\int_\Omega
U_{\theta,q_0}\phi^* dx}>0$. Due to
$\frac{a(x)\mu}{1+mu+k\mu}\rightarrow q_0(x)$ uniformly on
$\overline\Omega$ as $\mu\rightarrow\infty$, we know $U\rightarrow
U_{\theta,q_0}(x)$ uniformly on $\overline\Omega$. It follows from
(\ref{ee}) that
\begin{align}\eta&\le\int_{\Omega}[|\nabla\phi^*|^2-\theta{\phi^*}^2+2U{\phi^*}^2+\frac{a(x)\mu(1+k\mu)}{(1+mU+k\mu)^2}{\phi^*}^2]dx
\nonumber\\&=\eta^*+\int_\Omega[2(U-U_{\theta,q_0})+\frac{a(x)\mu(1+k\mu)}{(1+mu+k\mu)^2}-q_0(x)]{\phi^*}^2dx.
\end{align}
Thus $-\theta<\eta<M$ with $M>0$ independent of $\mu$. We claim that
$\lim\inf_{\mu\rightarrow\infty}\eta=r>0$. In fact, choose a
sequence $\mu_n\rightarrow\infty$ such that $\eta_n\rightarrow r$, and
\begin{eqnarray}\label{een}
-\Delta\phi_n=\theta\phi_n-2u_n\phi_n-\frac{a(x)
\mu_n(1+k\mu_n)}{(1+mu_n+k\mu_n)^2}\phi_n+\eta_n\phi_n~{\rm
in}~\Omega,~~\frac{\partial \phi_n}{\partial{n}}=0~{\rm
on}~\partial\Omega
\end{eqnarray} with normalized $\phi_n>0$, i.e. $\|\phi_n\|_2=1$. As $\int_\Omega|\nabla\phi_n|^2dx$ are
uniformly bounded with respect to $n$, there exists a subsequence
$\phi_{n_k}\rightharpoonup\phi_0$ weakly in $H^1(\Omega)$.
Obviously, $\phi_0\geq0$ and $\|\phi_0\|_2=1$. Multiply (\ref{een})
by $\varphi\in C_0^{\infty}(\Omega)$ and integrate by parts to have
$$\int_\Omega\nabla\phi_n\cdot\nabla\varphi dx=\int_\Omega[\theta\phi_n\varphi-2u_n\phi_n\varphi-\frac{a(x)\mu_n(1+k\mu_n)}{(1+mu_n+k\mu_n)^2}\phi_n\varphi+\eta_n\phi_n\varphi]dx.$$
Since $u_n\rightarrow U_{\theta,q_0}$ uniformly on
$\overline{\Omega}$ as $n\rightarrow\infty$, we have
$$\int_\Omega\nabla\phi_0\cdot\nabla\varphi
dx=\int_\Omega[\theta\phi_0\varphi-2U_{\theta,q_0}\phi_0\varphi-q_0(x)\phi_0\varphi+r\phi_0\varphi]dx.$$
Comparing with (\ref{eee}), we prove the claim that $r=\eta^*>0$.
So, there exists $\tilde{\mu}>0$ such that $\eta=\eta(\mu)>0$ when
$\mu>\tilde{\mu}$, which implies the linear stability of the
positive solutions to (\ref{su}).

Let
\begin{align}&H(t,u)=[MI-\Delta]^{-1}(M+\theta-u-t\frac{a(x)\mu}{1+mu+k\mu})u,\label{H}\\
&A=\{u\in C(\overline\Omega);~\varepsilon_0<u<\theta+1\}\nonumber
\end{align}
with $0<\varepsilon_0<\min_{x\in\overline{\Omega}}U_{\theta,q_0}(x)$, $0\leq
t\leq 1$, $M$ large. Define $$S(t,u)=u-H(t,u).$$It is easy to
see that $S(t,u)\neq0$ for all $u\in\partial A$, $0\leq t\leq 1$. For large $M$, by
the compactness of $H$, there are only finitely many
isolated fixed points in $A$, denoted by $u_1,\dots,u_m$. Together
with the linear stability of the positive solutions and the homotopy
invariance of fixed point index, we have
$$1=\mbox{index}(S(0,u),A,0)=\mbox{index}(S(1,u),A,0)=\sum_{i=1}^{m}\mbox{index}(H,u_i)=m.$$
Therefore, there is an unique positive fixed point to (\ref{H}) with
$t=1$ whenever $\mu>\tilde{\mu}$, i.e. problem (\ref{su}) has an unique positive solution. \qquad$\Box$\medskip

Now, we can deal with the uniqueness Theorem \ref{th3}.\medskip

\noindent {\bf Proof of Theorem \ref{th3}.} Let $(u,v)$ be a
positive solution of (\ref{s}) with large $\mu$. Linearize the
eigenvalue problem of (\ref{s}) at $(u,v)$ to have
 \begin{eqnarray}\label{31} \left\{
\begin{array}{llll}\displaystyle-\Delta\phi=\theta\phi-2u\phi-\frac{a(x)v(1+kv)}{(1+mu+kv)^2}\phi-\frac{a(x)u(1+mu)}{(1+mu+kv)^2}\psi
+\eta\phi&{\rm in}~\Omega,
\\ \displaystyle-\Delta\psi=\mu\psi-2v\psi+\frac{cu(1+mu)}{(1+mu+kv)^2}\psi+\frac{cv(1+kv)}{(1+mu+kv)^2}\phi+\eta\psi&{\rm
in}~\Omega_1,
\\
\displaystyle\frac{\partial\phi}{\partial n}=0~{\rm
on}~\partial\Omega,~~~~\frac{\partial\psi}{\partial n}=0~{\rm
on}~\partial\Omega_1.
\end{array}\right.
\end{eqnarray}
Here $\phi$, $\psi$ and $\eta$ may be complex-valued.

From Kato's inequality, we have
\begin{align}\label{32}
-\Delta|\phi|&\leq-{\rm Re}\Big(\frac{\bar{\phi}}{|\phi|}\Delta\phi\Big)\nonumber\\
&={\rm Re}\big(\theta|\phi|-2u|\phi|-\frac{a(x)v(1+kv)}{(1+mu+kv)^2}|\phi|+\frac{a(x)u(1+mu)}{(1+mu+kv)^2}\psi\cdot\frac{\bar{\phi}}{|\phi|}+\eta|\phi|\big)\nonumber\\
&\leq\theta|\phi|-2u|\phi|-\frac{a(x)v(1+kv)}{(1+mu+kv)^2}|\phi|+\frac{a(x)u(1+mu)}{(1+mu+kv)^2}|\psi|+{\rm Re}(\eta)|\phi|.
\end{align}
To obtain the linear stability, it suffices to prove that for any $\delta>0$, there
exists $\mu_{\delta}>0$ such that the eigenvalues $\eta$ of (\ref{31})
satisfy ${\rm Re}(\eta)\geq\eta^*-\delta$ when
$\mu\geq\mu_{\delta}$. Otherwise, there exist a $\delta_0>0$ and a
sequence $\{(\mu_n,\eta_n,u_n,v_n,\phi_n,\psi_n)\}_{n=1}^{\infty}$ satisfying
(\ref{31}) with $\|\phi_n\|_2+\|\psi_n\|_2=1$, and $\mu_n\rightarrow
\infty$ as $n\rightarrow\infty$ such that $\mbox{Re}(\eta_n)<\eta^*-\delta_0$. Replace $(\mu,\eta,u,v,\phi,\psi)$
in (\ref{32}) with $(\mu_n,\eta_n,u_n,v_n,\phi_n,\psi_n)$, multiply
by $|\phi_n|$, and then integrate by parts over $\Omega$ to have
\begin{align}\label{pp}
\int_\Omega|\nabla|\phi_n||^2dx\nonumber&\leq\int_\Omega\big(\theta|\phi_n|^2-2u_n|\phi_n|^2-\frac{a(x)v_n(1+kv_n)}{(1+mu_n+kv_n)}|\phi_n|^2\\
&~~+\frac{a(x)u_n(1+mu_n)}{(1+mu_n+kv_n)^2}|\psi_n||\phi_n|\big)
dx+(\eta^*-\delta_0)\int_\Omega|\phi_n|^2dx.
\end{align}
Let $r_n$ be the principal eigenvalue of the eigenvalue problem
\begin{align*}
-\Delta\varphi=\theta\varphi-2u_n\varphi-\frac{a(x)v_n(1+kv_n)}{(1+mu_n+kv_n)^2}\varphi+r_n\varphi~{\rm in}~\Omega,~~\frac{\partial\varphi}{\partial n}=0~{\rm on}~\partial\Omega.
\end{align*}
We know that
\begin{align*}r_n-\eta^*=\inf_{\varphi\in
H^1(\Omega)}\frac{\int_\Omega[|\nabla\varphi|^2-\theta\varphi^2+2u_n\varphi^2+\frac{a(x)v_n(1+kv_n)}{(1+mu_n+kv_n)^2}\varphi^2-\eta^*\varphi^2]dx}
{\int_\Omega\varphi^2dx},
\end{align*} and $r_n\rightarrow\eta^*$ by the proof of Proposition \ref{p3.1}. So, there exists a $N>0$ such that $r_n-\eta^*>-\frac{\delta_0}{2}$ for $n>N$. Thus by (\ref{pp}),
\begin{align*}
-\frac{\delta_0}{2}\int_\Omega|\phi_n|^2dx&<(r_n-\eta^*)\int_\Omega|\phi_n|^2dx
\nonumber\\&\leq\int_\Omega[|\nabla\phi_n|^2-\theta|\phi_n|^2+2u_n|\phi_n|^2+\frac{a(x)v_n(1+kv_n)}{(1+mu_n+kv_n)^2}|\phi_n|^2-\eta^*|\phi_n|^2]dx
\nonumber\\&\leq-\delta_0\int_{\Omega}|\phi_n|^2dx+\int_{\Omega_1}\frac{au_n(1+mu_n)}{(1+mu_n+kv_n)^2}|\psi_n||\phi_n|
dx.\end{align*}
Since $\frac{au_n(1+mu_n)}{(1+mu_n+kv_n)^2}\rightarrow 0$ in
$C(\overline\Omega_1)$ as $n\rightarrow \infty$, then
$\int_\Omega|\phi_n|^2dx\rightarrow 0$.

Using Kato's inequality again, we have
\begin{align*}-\Delta|\psi_n|\leq\mu_n|\psi_n|-2v_n|\psi_n|+\frac{cu_n(1+mu_n)}{(1+mu_n+kv_n)^2}|\psi_n|+
\frac{cv_n(1+kv_n)}{(1+mu_n+kv_n)^2}|\phi_n|+(\eta^*-\delta_0)|\psi_n|.
\end{align*}
Multiply by $|\psi_n|$ and integrate by parts over $\Omega_1$ to get
\begin{align*}\int_{\Omega_1}|\nabla|\psi_n||^2dx\nonumber&\leq\int_{\Omega_1}[\mu_n|\psi_n|^2-2v_n|\psi_n|^2+
\frac{cu_n(1+mu_n)}{(1+mu_n+kv_n)^2}|\psi_n|^2
\nonumber\\&~~+\frac{cv_n(1+kv_n)}{(1+mu_n+kv_n)^2}|\phi_n||\psi_n|+(\eta^*-\delta_0)|\psi_n|^2]dx
\nonumber\\&\leq\int_{\Omega_1}[-\mu_n+\frac{cu_n(1+mu_n)}{(1+mu_n+kv_n)^2}+\eta^*-\delta_0]|\psi_n|^2dx
\nonumber\\&~~+\int_{\Omega_1}\frac{cv_n(1+kv_n)}{(1+mu_n+kv_n)^2}|\phi_n||\psi_n|
dx.
\end{align*}
Consequently,
\begin{align*}\int_{\Omega_1}|\psi_n|^2dx\nonumber
&\leq\frac{1}{\mu_n}\int_{\Omega_1}[\frac{cu_n(1+mu_n)}{(1+mu_n+kv_n)^2}+\eta^*-\delta_0]|\psi_n|^2dx\\
&~~+\frac{1}{\mu_n}\int_{\Omega_1}\frac{cv_n(1+kv_n)}{(1+mu_n+kv_n)^2}|\phi_n||\psi_n|dx\nonumber\\
&\leq\frac{1}{\mu_n}\int_{\Omega_1}(\frac{c\theta}{1+m\theta}+\eta^*-\delta_0)|\psi_n|^2dx
+\frac{1}{\mu_n}\int_{\Omega_1}\frac{c}{k}|\phi_n||\psi_n|dx.
\end{align*}
This concludes $\int_{\Omega_1}|\psi_n|^2dx\rightarrow 0$ as
$n\rightarrow\infty$, since $\mu_n\rightarrow \infty$ and
$|\phi_n|,|\psi_n|$ are bounded in $L^2(\Omega_1)$.

In summary, we have obtained
$\int_\Omega|\phi_n|^2dx$,$\int_{\Omega_1}|\psi_n|^2dx\rightarrow
0$, as $n\rightarrow\infty$ which contradict with
$\|\phi_n\|_2+\|\psi_n\|_2=1$.

By using a similar argument to that in the proof of Proposition
\ref{p3.1}, we get from the linear stability of the positive
solutions to (\ref{s}) and  Proposition \ref{p3.1} that the
solution of (\ref{s}) must be unique when $\mu>\max\{\tilde{\mu},\mu_0\}$ with
$\mu_0=\inf\{\mu_\delta; \delta\in (0,\eta^*)\}$.

Finally, we consider the asymptotic behavior of the unique positive
solution $(u,v)$ as $\mu\rightarrow \infty$. Since
$\frac{cu}{1+mu+kv}\le \frac{c\theta}{1+m\theta+k\mu}\rightarrow
0$ as $\mu\rightarrow \infty$, for any $\epsilon>0$, there is a
$\mu_\epsilon>0$ such that $\frac{cu}{1+mu+kv}<\epsilon$ for
$\mu>\mu_\epsilon$. Then
\begin{align*}
\mu v-v^2\leq-\Delta v\leq(\mu+\epsilon)v-v^2~{\rm in}~\Omega_1,
\end{align*} which yields $\mu\leq
v\leq\mu+\epsilon$ for $\mu>\mu_\epsilon$. And thus $v-\mu\rightarrow 0$
as $\mu\rightarrow\infty$. We know that
$\frac{a(x)v}{1+mu+kv}\rightarrow q_0(x)$, and then $u\rightarrow
U_{\theta,q_0}(x)$ uniformly on $\overline\Omega$ as
$\mu\rightarrow\infty$. \qquad$\Box$
\section{Discussion}
\mbox\indent
In a reaction-diffusion system of predator-prey PDE model,
in addition to the interaction mechanism between the species,  the
behavior of the species is also affected by the diffusion of the
species, as well as the size and geometry of the habitat. Obviously,
the prey species would die out under excessive predation from nature
or humans. The results obtained in this paper show the way in which the
created protection zone saves the endangered prey species in the
diffusive predator-prey model with Beddington-DeAngelis type functional
response and non-flux boundary conditions.

Compared with previous results on protection zone problems with other functional
responses such as the Lotka-Voltera type
competition system \cite{DL}, Holling II type predator-prey system
\cite{RC,DS2006}, and Leslie type predator-prey system \cite{DPW}, richer dynamic properties have been  observed for the model (\ref{a}) with B-D type functional
response in this paper. It can be found that a total of four threshold values are obtained here
for the prey birth rate $\theta$, i.e., $\theta_0$, $\theta^*$, $\theta_1$ (for
the predator growth rate $\mu>0$) and $\theta_*$ (for $\mu\le 0$).

The first threshold value $\theta_0$ gives the necessary condition
for establishing a protection zone to save the prey $u$.  By Theorem
\ref{th1}, the survival of $u$ could be automatically ensured
without protection zones whenever
$\theta>\theta_0=\frac{a}{k}$, which can be realized when the
refuge ability of the prey is properly large that
$k>\frac{a}{\theta}$, or the predation rate is small that
$a<{\theta}{k}$. In other words, the protection zones have to be
made only if the prey's refuge ability is too weak with respect to its
birth rate $\theta$ and the predation rate $a$. This matches with
the mechanistic derivation for the B-D type functional response
proposed in \cite{GG}. In addition, Theorem \ref{th1} says also that
the protection zones are unnecessary if the predator's growth rate
$\mu\le 0$, where the predator species $v$ can not live without the
prey $u$, and thus the extinction of $v$ cannot take place after of $u$.

The second threshold value is
 $\theta^*=\theta^*(\mu,\Omega_0)=\lambda_1(q(x))$ with
$q(x)=\frac{a(x)\mu}{1+k\mu}$ and $\mu>0$. By Theorem \ref{th2},
the positive steady states can be attained  for
$\theta\in(\theta^*,\theta_0)$. Due to the monotonicity of the
principal eigenvalue $\lambda_1=\lambda_1(q(x))$ with respect to
$q(x)$, we know that the threshold value $\theta^*$ would be
enlarged (and hence harmful for the prey $u$) when the predation
rate $a(x)$ or the predator's growth rate $\mu$ increase, or when the
prey refuge $k$ or the size of the protection zone $\Omega_0$
decrease. Conversely, Theorem \ref{th2} says also that the prey
$u$ must become extinct when $\theta\le\theta^*$ with the
handling time $m$ of $v_H$ being shorter than $\frac{(k\mu+1)^2}{a\mu}$.
All of these match with those in \cite{GG}. In addition, since $\theta^*{(\mu,\Omega_0)}\le\theta_0$ is strictly increasing with respect to $\mu$ and decreasing when enlarging $\Omega_0$, letting $\mu\rightarrow\infty$, we get the third threshold value $\theta_1=\theta_1(\Omega_0)$  such that if $\theta>\theta_1(\Omega_0)$, prey species could survive no matter how large the predator's growth rate is. The critical $\theta=\theta_1(\Omega_0)$ implies a critical size of the protection zone as well, namely, if the real protection zone $\widetilde{\Omega}_0\Supset \Omega_0$, the survival of the prey with such birth rate $\theta$ is independent of the predator's growth rate.  Also, the uniqueness and linear stability obtained in Theorem \ref{th3} for $\mu$ large enough are reasonable because $\frac{cu}{1+mu+kv}\rightarrow 0$, and hence $v-\mu\rightarrow 0$ as $\mu\rightarrow\infty$.

Since the condition $\mu \le 0$ yields the survival of $u$ without protection
zones by Theorem \ref{th1}, the fourth threshold value $\theta_*$ obtained in Theorem \ref{th2'} with $\mu\in (-\frac{c}{m},0]$ is only made for $v$ alive. In fact, the conversion of prey is limited by $\frac{c}{m}$, as shown in the proof of Theorem \ref{th1}, the predator $v$ must be die out if its growth rate $\mu\le-\frac{c}{m}$.
With such non-positive growth rate $\mu\in (-\frac{c}{m},0]$, there should be properly large number of prey to survive the predator, just as described via the criterion
$\theta>\theta_*=\frac{|\mu|}{c-m|\mu|}\ge 0$ in Theorem \ref{th2'}. It is worth pointing out that the threshold value $\theta_*$ for alive
predator $v$ would be enlarged as $m$ (the handling time of $v_H$)
is increasing. This well matches the mechanism in  \cite{GG}.

We have shown the effect of refuge ability of the
prey and protection zones on the the coexistence and stability of
predator-prey species in the paper. In fact, protection zone can be
regarded as another refuge offered by human intervention, which is necessary if the prey's refuge ability is
too weak in the predator-prey system to prevent the extinction of the prey populations.  The critical sizes of protection zones, obtained in this paper and
represented by the principal eigenvalues $\theta^*=\lambda_1(q(x))$ and
$\theta_1=\lambda_1(q_0(x))$, show the basic requirement (depending on the predator's growth rate) and the sufficient one (working under any predator's growth rate), respectively. The results of the present paper would be helpful to design
nature reserves and no-fishing zones, etc.



\end{document}